# ASYMPTOTIC NORMALITY OF EXTREME
# VALUE ESTIMATORS ON $C[0,1]$

BY JOHN H. J. EINMAHL AND TAO LIN

*Tilburg University and Xiamen University*

Consider $n$ i.i.d. random elements on $C[0,1]$. We show that, under an appropriate strengthening of the domain of attraction condition, natural estimators of the extreme-value index, which is now a continuous function, and the normalizing functions have a Gaussian process as limiting distribution. A key tool is the weak convergence of a weighted tail empirical process, which makes it possible to obtain the results uniformly on $[0,1]$. Detailed examples are also presented.

**1. Introduction.** Recently considerable progress has been made in the interesting field of infinite-dimensional extreme value theory, where the data are (continuous) functions. After the characterization of max-stable stochastic processes in $C[0,1]$ by Giné, Hahn and Vatan [11], de Haan and Lin [4, 5] investigated the domain of attraction conditions and established weak consistency of estimators of the extreme value index, the centering and standardizing sequences, and the exponent measure.

Statistics of infinite-dimensional extremes will find various applications, for example, to coast protection (flooding) and risk assessment in finance. For an application to coast protection, consider the northern part of The Netherlands, which lies for a substantial part below sea level. Since there is no natural coast defense there, the area is protected by a long dike against inundation. Flooding of the dike at any place could lead to flooding of the whole area, so the approach via function-valued data is the appropriate one here. In finance, the intra-day return of a stock is defined as the ratio of the price of a stock at a certain time $t$ during the day to the price at market opening. This process can be well described, when we measure time in days, with a continuous function on $[0,1]$. For various risk analysis problems (e.g., problems dealing with options), intra-day returns of the stock need to be taken into account, instead of just the daily returns (i.e., the function values at 1). Sampling on $n$ days puts us in a position to apply statistics of extremes to these problems.

Also, from a mathematical point of view, the research is challenging, because of the new features of $C[0,1]$-valued random elements, when compared to random variables or vectors, in particular, the uniformity in $t \in [0,1]$ of the results asks for novel approaches.

It is the purpose of this paper to establish asymptotic normality of estimators of the extreme value index, which is now an element of $C[0,1]$, and of the normalizing sequences. In fact, we will show the asymptotic normality on $C[0,1]$ of the estimators under a suitable second-order condition and present all the limiting processes involved in terms of one underlying Wiener process, which means that we have the simultaneous weak convergence of all the estimators. The results are, on the one hand, interesting in themselves, because the extreme value index measures the tail heaviness of the distribution of the data, and, on the other hand, the results are a major step toward the estimation of probabilities of rare events in $C[0,1]$; see [7] for a study of this problem in the finite-dimensional case.

In order to be more explicit let us now specify the setup and introduce notation. Let $\xi_1, \xi_2, \ldots$ be i.i.d. random elements on $C[0,1]$. Define $F_t : \mathbb{R} \to [0,1]$ by $F_t(x) = P\{\xi_i(t) \leq x\}$. Throughout assume that

$$P\left\{ \inf_{t \in [0,1]} \xi_i(t) > 0 \right\} = 1 \tag{1}$$

and that

(2)     $F_t$ is a continuous and strictly increasing function on its support.

Define

$$U_t(s) = F_t^{\leftarrow}(1 - 1/s), \qquad s > 0, 0 \leq t \leq 1.$$

We assume that the domain of attraction condition holds, that is,

$$\left\{ \left( \max_{i=1,\ldots,n} \xi_i(t) - b_t(n) \right) \Big/ a_t(n), t \in [0,1] \right\} \text{ converges in distribution} \tag{3}$$

on $C[0,1]$ to a stochastic process, $\eta$, say, with nondegenerate marginals, where $a_t(n) > 0$ and $b_t(n)$ are continuous (in $t$) normalizing functions, chosen in such a way that, for each $t$,

$$P\{\eta(t) \leq x\} = \exp(-(1 + \gamma(t)x)^{-1/\gamma(t)});$$

see [4]. We can and will take $b_t \equiv U_t$. Then $\gamma:[0,1] \to \mathbb{R}$, the extreme value index (function), is continuous. Define

$$\zeta_i(t) = \frac{1}{1 - F_t(\xi_i(t))},$$

$$\bar{\eta}(t) = (1 + \gamma(t)\eta(t))^{1/\gamma(t)},$$

and $\nu_s(E) = sP\{\zeta_i \in sE\}$, with $sE = \{sh : h \in E\}$. Clearly the $\zeta_i(t)$ are standard Pareto random variables, that is, $P\{\zeta_i(t) \leq x\} = 1 - 1/x$, $x \geq 1$. It follows from Theorem 2.8 in [4] that there exists a measure $\nu$ on $C[0,1]$ that is homogeneous [i.e., for a Borel set $A$ and $r > 0$, $\nu(rA) = \frac{1}{r}\nu(A)$], such that for any positive $g \in C[0,1]$ and compact set $K \subset [0,1]$,

$$P\{\bar{\eta}(t) < g(t), \text{ for all } t \in K\}$$
$$= \exp(-\nu(\{f \in C[0,1], f(t) \geq g(t), \text{ for some } t \in K\}))$$

and

(4) $$\nu_s \to \nu \qquad \text{as } s \to \infty,$$

weakly on $S_c := \{f \in C[0,1] : \sup_{t \in [0,1]} f(t) \geq c\}$ for any $c > 0$, and

$$\frac{U_t(nx) - U_t(n)}{a_t(n)} \to \frac{x^{\gamma(t)} - 1}{\gamma(t)} \qquad \text{as } n \to \infty,$$

uniformly in $t \in [0,1]$ and locally uniformly in $x \in (0, \infty)$.

*Throughout* we assume that $k = k(n) \in \{1, \ldots, n-1\}$ is a sequence of positive integers satisfying

(5) $$k \to \infty \quad \text{and} \quad \frac{k}{n} \to 0 \qquad \text{as } n \to \infty.$$

Fix $t \in [0,1]$. Let $\xi_{1,n}(t) \leq \xi_{2,n}(t) \leq \cdots \leq \xi_{n,n}(t)$ be the order statistics of $\xi_i(t), i = 1, 2, \ldots, n$. We define the following statistical functions:

(6) $$M_n^{(r)}(t) = \frac{1}{k}\sum_{i=0}^{k-1}(\log \xi_{n-i,n}(t) - \log \xi_{n-k,n}(t))^r, \qquad r = 1, 2.$$

Set $\gamma^+(t) = \gamma(t) \vee 0$, $\gamma^-(t) = \gamma(t) \wedge 0$ and observe that $\gamma(t) = \gamma^+(t) + \gamma^-(t)$. Now we define estimators for $\gamma^+(t)$, $\gamma^-(t)$, $\gamma(t)$, $a_t(\frac{n}{k})$ and $b_t(\frac{n}{k})$ as in [9]:

(7) $$\hat{\gamma}_n^+(t) = M_n^{(1)}(t) \qquad \text{(Hill estimator)};$$

(8) $$\hat{\gamma}_n^-(t) = 1 - \frac{1}{2}\left(1 - \frac{(M_n^{(1)}(t))^2}{M_n^{(2)}(t)}\right)^{-1};$$

(9) $$\hat{\gamma}_n(t) = \hat{\gamma}_n^+(t) + \hat{\gamma}_n^-(t) \qquad \text{(moment estimator)};$$

(10) $\quad \hat{U}_t\left(\dfrac{n}{k}\right) = \xi_{n-k,n}(t) \qquad\qquad\qquad$ (location estimator);

(11) $\quad \hat{a}_t\left(\dfrac{n}{k}\right) = \xi_{n-k,n}(t)\hat{\gamma}_n^+(t)(1-\hat{\gamma}_n^-(t)) \qquad$ (scale estimator).

For fixed $t$ these are well-known one-dimensional estimators. Observe that $\hat{\gamma}_n^+(t)$ and $\hat{\gamma}_n^-(t)$ are not equal to $(\hat{\gamma}_n)^+(t)$ and $(\hat{\gamma}_n)^-(t)$, respectively.

The following weak consistency results have been shown in [5].

THEOREM 1.1. *If* (1), (2), (3) *and* (5) *hold, then*

(12) $$\sup_{0\le t\le 1} |\hat{\gamma}_n^+(t) - \gamma^+(t)| \xrightarrow{P} 0,$$

(13) $$\sup_{0\le t\le 1} |\hat{\gamma}_n(t) - \gamma(t)| \xrightarrow{P} 0,$$

(14) $$\sup_{0\le t\le 1} \left|\frac{\hat{U}_t(n/k) - U_t(n/k)}{a_t(n/k)}\right| \xrightarrow{P} 0,$$

(15) $$\sup_{0\le t\le 1} \left|\frac{\hat{a}_t(n/k)}{a_t(n/k)} - 1\right| \xrightarrow{P} 0.$$

The main results of the paper and examples are presented in Section 2; the proofs are deferred to Section 3.

**2. Main results.** In this section we present our main result, dealing with the asymptotic normality of the estimators of which the weak consistency is shown in Theorem 1.1. In order to establish our main result, we first present a result that is a key tool for its proof. This result deals with the weak convergence of a tail empirical process based on the $\zeta_i, i = 1,\ldots,n$. Write $C_{t,x} = \{h \in C[0,1] : h(t) \ge x\}$ and define

$$S_{n,t}(x) = \frac{1}{n}\sum_{i=1}^n \mathbb{1}_{\{\zeta_i \in C_{t,x}\}} = \frac{1}{n}\sum_{i=1}^n \mathbb{1}_{\{\zeta_i(t) \ge x\}}.$$

Denote, with $k$ as in (5), the corresponding tail empirical process with

$$w_n(t,x) = \sqrt{k}\left(\frac{n}{k}S_{n,t}\left(x\frac{n}{k}\right) - \frac{1}{x}\right).$$

Let $c > 0$ and define $\mathcal{C} = \{C_{t,x} : 0 \le t \le 1, x \ge c\}$. Let $W$ be a zero-mean Gaussian process defined on $\mathcal{C}$ with $EW(C_{t,x})W(C_{s,y}) = \nu(C_{t,x} \cap C_{s,y})$. Clearly, for fixed $t \in [0,1]$, $\{W(C_{t,1/y}), y \le \frac{1}{c}\}$ is a standard Wiener process, since $\nu(C_{t,1/y_1} \cap C_{t,1/y_2}) = \nu(C_{t,1/(y_1 \wedge y_2)}) = y_1 \wedge y_2$. For $\beta > 0$, set, for any $(t,x)$,

$(s, y) \in [0, 1] \times [c, \infty)$,

$$d((t,x),(s,y))$$
$$= \sqrt{E(x^\beta W(C_{t,x}) - y^\beta W(C_{s,y}))^2}$$
$$= \sqrt{(x^\beta - y^\beta)^2 \nu(C_{t,x} \cap C_{s,y}) + x^{2\beta}\nu(C_{t,x} \setminus C_{s,y}) + y^{2\beta}\nu(C_{s,y} \setminus C_{t,x})}.$$

Observe that (4) implies that $\frac{n}{k}P\{\zeta_i(t) \geq \frac{n}{k}x, \zeta_i(s) \geq \frac{n}{k}y\} = \nu_{n/k}(C_{t,x} \cap C_{s,y}) \to \nu(C_{t,x} \cap C_{s,y})$.

Define, for $K > 0$,

$$E_{s,\delta} = \left\{ h \in C[0,1] : h \geq 0, \frac{|h(t) - h(s)|}{h(s)} \leq K\left(\log\frac{1}{\delta}\right)^{-3} \text{ for all } t \in [s, s+\delta] \right\}$$

and assume that, for $0 \leq \beta < \frac{1}{2}$ and all $c_1 > 0$, for some $K > 0$ and for large enough $v$ there exists a $\delta_0 > 0$ such that for all $\delta \in (0, \delta_0]$,

$$(16) \quad \sup_{s \in [0,1]} P\left\{ \zeta_1 \notin E_{s,\delta} \,\Big|\, \sup_{t \in [s, s+\delta]} \zeta_1(t) \geq v \right\} \leq c_1 \left(\log\frac{1}{\delta}\right)^{-(2+2\beta)/(1-2\beta)}.$$

For convenient presentation and convenient application in the proofs of the main result, this result is presented in an approximation setting (with the random elements involved defined on *one* probability space), via the Skorohod–Dudley–Wichura construction. So the random elements $w_n$ and $W$ in this theorem are only equal in distribution to the original ones, but we do not add the usual tildes to the notation.

THEOREM 2.1. *Let $0 \leq \beta < \frac{1}{2}$. Suppose the conditions (2), (3), (5) and (16) hold. Then for the new $w_n$ and $W$ mentioned above, we have, for any $c > 0$,*

$$(17) \quad \sup_{0 \leq t \leq 1, x \geq c} x^\beta |w_n(t,x) - W(C_{t,x})| \xrightarrow{P} 0 \quad \text{as } n \to \infty.$$

*Define $Z(t,x) = x^\beta W(C_{t,x})$. Then the process $Z$ is bounded and uniformly $d$-continuous on $[0,1] \times [c, \infty)$.*

Note that it is well known that, for *one* fixed $t$, the restriction $\beta < \frac{1}{2}$ is also necessary for weak convergence of the (one-dimensional) tail empirical process. So our condition on $\beta$ in the present infinite-dimensional setting is the same as in dimension one.

Condition (16) is needed to prove tightness. It prevents the continuous random function from having extremely large oscillations. From the examples below we see that it is a rather weak condition, since they amply satisfy this condition.

It is important to transform the $\xi_i$ to processes with standard marginals, as we did by transforming to the $\zeta_i$. Although the choice of standard Pareto marginals is convenient, it is also reasonable to transform to other marginal distributions, such as the uniform-$(0,1)$ distribution. Clearly, uniform-$(0,1)$ marginals are obtained by taking $1/\zeta_i$. It is interesting to note and readily checked that the set $E_{s,\delta}$ used in condition (16) is invariant under this transformation.

It can be useful to replace $|h(t) - h(s)|/h(s)$ by $|\log h(t) - \log h(s)|$ in the definition of $E_{s,\delta}$ used in condition (16). The thus obtained version of condition (16) is equivalent to the stated one ($K$ can be different), but it might be easier to check for certain processes.

We also need the following corollary, which deals with certain quantiles and can be obtained by the usual "inversion" from the tail empirical process theorem.

COROLLARY 2.2. *Let $\alpha \in C[0,1]$. We have, under the conditions of Theorem* 2.1,

$$\text{(18)} \quad \sup_{0 \leq t \leq 1} \left| \sqrt{k} \left( \left( \zeta_{n-k,n}(t) \frac{k}{n} \right)^{\alpha(t)} - 1 \right) - \alpha(t) W(C_{t,1}) \right| \xrightarrow{P} 0 \quad \text{as } n \to \infty.$$

Finally, we present the main result, which gives the asymptotic distributions of the estimators of $\gamma^+$, $\gamma$, $a$ and $b$ in terms of the process $W$, figuring in Theorem 2.1.

THEOREM 2.3. *Suppose the conditions of Theorem* 2.1 *and* (1) *are satisfied and the following second-order condition holds: for every $t \in [0,1]$, there exists a function $A_t$ not changing sign near infinity with $\lim_{v \to \infty} \sup_{0 \leq t \leq 1} |A_t(v)| = 0$, such that, as $v \to \infty$,*

$$\text{(19)} \quad \left( \frac{\log U_t(vx) - \log U_t(v)}{a_t(v)/U_t(v)} - \frac{x^{\gamma^-(t)} - 1}{\gamma^-(t)} \right) \Big/ A_t(v) \to H_{\gamma^-(t), \rho(t)}(x),$$

*uniformly in $t \in [0,1]$ and locally uniformly in $x > 0$, with*

$$H_{\gamma^-(t),\rho(t)}(x) = \int_1^x y^{\gamma^-(t)-1} \int_1^y u^{\rho(t)-1}\, du\, dy,$$

*with $\rho(t) \in [-\infty, 0]$ for all $t \in [0,1]$.*

*If, as $n \to \infty$,*

$$\text{(20)} \quad \sqrt{k} \sup_{0 \leq t \leq 1} \left| A_t\!\left(\frac{n}{k}\right) \right| \to 0$$

*and*

$$\text{(21)} \quad \sqrt{k} \sup_{0 \leq t \leq 1} \left| \frac{a_t(n/k)}{U_t(n/k)} - \gamma^+(t) \right| \to 0,$$

*then we have*

$$\sup_{0\leq t\leq 1} |\sqrt{k}(\hat{\gamma}_n^+(t) - \gamma^+(t)) - \gamma^+(t)\mathcal{P}(t)| \xrightarrow{P} 0, \tag{22}$$

$$\sup_{0\leq t\leq 1} |\sqrt{k}(\hat{\gamma}_n(t) - \gamma(t)) - \Gamma(t)| \xrightarrow{P} 0, \tag{23}$$

$$\sup_{0\leq t\leq 1} \left|\sqrt{k}\frac{\hat{U}_t(n/k) - U_t(n/k)}{a_t(n/k)} - \mathcal{U}(t)\right| \xrightarrow{P} 0, \tag{24}$$

$$\sup_{0\leq t\leq 1} \left|\sqrt{k}\left(\frac{\hat{a}_t(n/k)}{a_t(n/k)} - 1\right) - \mathcal{A}(t)\right| \xrightarrow{P} 0, \tag{25}$$

*where $\mathcal{P}, \Gamma, \mathcal{U}$ and $\mathcal{A}$ are defined in terms of the process $W$ as follows:*

$$\mathcal{P}(t) = \int_1^\infty W(C_{t,x}) \frac{dx}{x^{1-\gamma^-(t)}} - \frac{1}{1-\gamma^-(t)} W(C_{t,1}),$$

$$\mathcal{Q}(t) = 2\int_1^\infty W(C_{t,x}) \frac{x^{\gamma^-(t)} - 1}{\gamma^-(t)} \frac{dx}{x^{1-\gamma^-(t)}}$$
$$\quad - 2((1-\gamma^-(t))(1-2\gamma^-(t)))^{-1} W(C_{t,1}),$$

$$\Gamma(t) = \{\gamma^+(t) - 2(1-\gamma^-(t))^2(1-2\gamma^-(t))\}\mathcal{P}(t)$$
$$\quad + \tfrac{1}{2}(1-\gamma^-(t))^2(1-2\gamma^-(t))^2 \mathcal{Q}(t),$$

$$\mathcal{U}(t) = W(C_{t,1}),$$

$$\mathcal{A}(t) = \gamma(t)W(C_{t,1}) + (3-4\gamma^-(t))(1-\gamma^-(t))\mathcal{P}(t)$$
$$\quad - \tfrac{1}{2}(1-\gamma^-(t))(1-2\gamma^-(t))^2 \mathcal{Q}(t), \qquad t \in [0,1].$$

Condition (19) is a uniform version of one of the natural, well-studied second-order conditions of univariate extreme value theory; see [3] and [8]. For $\rho(t) > -\infty$, the absolute value of the function $A_t$ is regularly varying of order $\rho(t)$ and specifies the rate of convergence in (3). For $\rho(t) = -\infty$, we can choose $A_t$ such that its absolute value tends to 0 faster than a given power function. Large values of $|\rho(t)|$ yield fast convergence, whereas small values and, in particular, the case $\rho(t) = 0$, correspond to (very) slow convergence.

Note that for the case $\inf_{t\in[0,1]} \gamma(t) > 0$ and $\sup_{t\in[0,1]} \rho(t) < 0$, it follows from the second-order condition (19) that

$$\sup_{t\in[0,1]} \left|\frac{\rho(t)((a_t(v)/U_t(v)) - \gamma^+(t))}{A_t(v)} - 1\right| \to 0 \qquad \text{as } v \to \infty. \tag{26}$$

So in this case (21) is superfluous, since it follows from (26) and (20). Also, note that condition (20) can be replaced by the stronger, but easier-to-check condition: for some $\varepsilon > 0$,

$$\sqrt{k}\left(\frac{n}{k}\right)^{\varepsilon + \sup_{t \in [0,1]} \rho(t)} \to 0.$$

For the case $\sup_{t \in [0,1]} \gamma(t) < 0$ and $\sup_{t \in [0,1]} \rho(t) < 0$, it follows from the second-order condition (19) that conditions (20) and (21) can be replaced by the stronger condition: for some $\varepsilon > 0$,

$$\sqrt{k}\left(\frac{n}{k}\right)^{\varepsilon + \sup_{t \in [0,1]} \rho(t) \vee \sup_{t \in [0,1]} \gamma(t)} \to 0.$$

When $\sup_{t \in [0,1]} \rho(t) = 0$ or $\gamma(t_1) = 0$ for some $t_1 \in [0, 1]$ [this also implies $\rho(t_1) = 0$], we do not have a simple sufficient condition on the growth of $k$, but it is necessary that $k$ grow more slowly than any power of $n$.

EXAMPLES. In order to illustrate the theory and for a better understanding of the conditions in the theorem, we consider two classes of examples.

Let $f$ be a unimodal, continuous probability density function on the real line satisfying, for some $K, \delta_0 > 0$ and for all $\delta \leq \delta_0$,

$$\sup_{|t-s| \leq \delta} \frac{|f(t) - f(s)|}{f(s)} \leq K\left(\log \frac{1}{\delta}\right)^{-3}.$$

This condition is satisfied by, for example, the double exponential density $((\lambda/2)\exp(-\lambda|x|))$ or the $t$-distribution for any number of degrees of freedom. Let $(X_j, Y_j), j = 1, 2, \ldots$, be an enumeration of the points of a homogeneous, rate 1, Poisson process on $\mathbb{R} \times \mathbb{R}^+$. Now define

$$\xi_1(t) = \sup_j \frac{f(t + X_j)}{Y_j} \qquad \text{for } t \in [0, 1].$$

This process is studied in detail in [1]; see also [6]. In particular, $\xi_1$ is a continuous, stationary, max-stable (i.e., limiting) process with marginals $F_t(x) = P\{\xi_1(t) \leq x\} = \exp(-1/x), x > 0$ (for all $t$ and $i$). Observe that $\gamma \equiv 1$ here.

For this process we will only check condition (16) in detail. The other conditions are easily seen to hold. In particular, (3) holds since the process is max-stable and (19) is well known to hold for the distribution function $F_t$ in the univariate case, that is, for fixed $t$. Since $F_t$ does not depend on $t$, it therefore holds uniformly in $t$ and, hence, (19) holds. Note that $\rho \equiv -1$. We now check (16). We have that

$$\frac{1}{1 - F_t(x)} = \frac{1}{1 - \exp(-1/x)} =: g(x).$$

Hence, for large values of $x$ the transformation $g$ is close to the identity. Note that $g(x) \geq x$ and $g'(x) \leq 1$ for $x > 0$. Hence, for $s, t$ and $\delta$ as above,

$$\frac{|g(\xi_1(t)) - g(\xi_1(s))|}{g(\xi_1(s))} \leq \frac{|g(\sup_j f(t+X_j)/Y_j) - g(\sup_j f(s+X_j)/Y_j)|}{g(\sup_j f(s+X_j)/Y_j)}$$

$$\leq \frac{|\sup_j f(t+X_j)/Y_j - \sup_j f(s+X_j)/Y_j|}{\sup_j f(s+X_j)/Y_j}$$

$$\leq \frac{\sup_j |f(t+X_j)/Y_j - f(s+X_j)/Y_j|}{\sup_j f(s+X_j)/Y_j}$$

$$\leq \sup_j \frac{|f(t+X_j)/Y_j - f(s+X_j)/Y_j|}{f(s+X_j)/Y_j}$$

$$\leq \sup_j \frac{|f(t+X_j) - f(s+X_j)|}{f(s+X_j)}$$

$$\leq K \left(\log \frac{1}{\delta}\right)^{-3}.$$

So condition (16) is satisfied since the probability involved in the condition is equal to 0 for all $s \in [0,1]$.

Let $Y$ be a standard Pareto random variable, that is, $P\{Y \leq x\} = 1 - 1/x$ for $x \geq 1$, and let $B$ be a random element of $C[0,1]$ such that $B(t) > 0$, $EB(t) = 1$ for all $t \in [0,1]$, and $E\sup_t B(t) < \infty$. Assume $Y$ and $B$ are independent. Define

$$\xi_1(t) = YB(t) \qquad \text{for } t \in [0,1];$$

see also [12]. We first show that $\xi_1$ satisfies the domain of attraction condition (3); more precisely, we show that $\frac{1}{n} \max_{i=1,\ldots,n} \xi_i$ converges in distribution to $\eta$, where $\xi_1, \ldots, \xi_n$ are i.i.d. We need to show the convergence of the finite-dimensional distributions and tightness. For the convergence of the finite-dimensional distributions, let $t_1, \ldots, t_k \in [0,1]$, $x_1, \ldots, x_k \geq 0$ and $\max_{j=1,\ldots,k} x_j > 0$. Now we have

$$\log P\left\{\frac{1}{n} \max_{i=1,\ldots,n} Y_i B_i(t_1) \leq x_1, \ldots, \frac{1}{n} \max_{i=1,\ldots,n} Y_i B_i(t_k) \leq x_k\right\}$$

$$= n \log P\left\{\frac{1}{n} YB(t_1) \leq x_1, \ldots, \frac{1}{n} YB(t_k) \leq x_k\right\}$$

$$= n \log \left[1 - P\left\{\frac{1}{n} YB(t_1) > x_1 \text{ or } \ldots \text{ or } \frac{1}{n} YB(t_k) > x_k\right\}\right]$$

$$\sim -nP\left\{\frac{1}{n} YB(t_1) > x_1 \text{ or } \ldots \text{ or } \frac{1}{n} YB(t_k) > x_k\right\}$$

$$= -nP\left\{Y > \min_{j=1,\ldots,k} \frac{nx_j}{B(t_j)}\right\}$$

$$= -nE\left\{\left(\max_{j=1,\ldots,k} \frac{B(t_j)}{nx_j}\right) \wedge 1\right\}$$

$$= -E\left\{\left(\max_{j=1,\ldots,k} \frac{B(t_j)}{x_j}\right) \wedge n\right\}$$

$$\to -E\left\{\max_{j=1,\ldots,k} \frac{B(t_j)}{x_j}\right\} \quad \text{as } n \to \infty.$$

This settles the convergence of the finite-dimensional distributions. Note that for $k=1$ the last expression is simply $-1/x_1$, which means that again $\gamma \equiv 1$. Next we consider the tightness. From the derivation above, it follows that $P\{\frac{1}{n}\max_{i=1,\ldots,n} Y_i B_i(0) > M\}$ can be made arbitrarily small for $M$ and $n$ large enough. So it remains to show that for $\varepsilon > 0$ there exists a $\delta > 0$ such that, for large $n$ large enough,

$$P\left\{\sup_{|t-s|<\delta}\left|\frac{1}{n}\max_{i=1,\ldots,n} Y_i B_i(t) - \frac{1}{n}\max_{i=1,\ldots,n} Y_i B_i(s)\right| > \varepsilon\right\} < \varepsilon.$$

We have

$$\frac{1}{\varepsilon}P\left\{\sup_{|t-s|<\delta}\left|\frac{1}{n}\max_{i=1,\ldots,n} Y_i B_i(t) - \frac{1}{n}\max_{i=1,\ldots,n} Y_i B_i(s)\right| > \varepsilon\right\}$$

$$\leq \frac{1}{\varepsilon}P\left\{\max_{i=1,\ldots,n}\sup_{|t-s|<\delta}\left|\frac{1}{n}Y_i B_i(t) - \frac{1}{n}Y_i B_i(s)\right| > \varepsilon\right\}$$

$$\leq \frac{n}{\varepsilon}P\left\{\sup_{|t-s|<\delta}\left|\frac{1}{n}YB(t) - \frac{1}{n}YB(s)\right| > \varepsilon\right\}$$

$$= \frac{n}{\varepsilon}EP\left\{Y > \frac{n\varepsilon}{\sup_{|t-s|<\delta}|B(t)-B(s)|}\bigg|\sup_{|t-s|<\delta}|B(t)-B(s)|\right\}$$

$$\leq \frac{n}{\varepsilon}E\left\{\frac{\sup_{|t-s|<\delta}|B(t)-B(s)|}{n\varepsilon}\right\}$$

$$= \frac{1}{\varepsilon^2}E\left\{\sup_{|t-s|<\delta}|B(t)-B(s)|\right\}.$$

Since $B \in C[0,1]$, we have $\sup_{|t-s|<\delta}|B(t)-B(s)| \to 0$. But $\sup_{|t-s|<\delta}|B(t)-B(s)| \leq 2\sup_t B(t)$ and, by assumption, $E\sup_t B(t) < \infty$, so by Lebesgue's dominated convergence theorem $E\{\sup_{|t-s|<\delta}|B(t)-B(s)|\} \to 0$ as $\delta \downarrow 0$. This completes the proof of the tightness.

In the sequel we will make the specific choice $B(t) = \exp(W(t) - \frac{t}{2}), t \in [0,1]$, with $W$ a standard Wiener process. $B$ is a geometric Brownian motion. Note that this process satisfies the conditions on $B$ specified in the beginning of this example. In particular, $E \sup_t B(t) < \infty$ follows from simply bounding $W(t) - \frac{t}{2}$ by $W(t)$ and the fact that the distribution function of $\sup_t W(t)$ is well known to be $2\Phi - 1$, where $\Phi$ is the standard normal distribution function. The corresponding process $\xi_1 = YB$ has been introduced in [12]. It remains to consider (19) and (16) for this process. It follows from a straightforward calculation that, for every $M > 1$, uniformly in $t \in [0,1]$,

$$\frac{1}{u} \geq P\left\{Y \exp\left(W(t) - \frac{t}{2}\right) > u\right\} \geq \frac{1}{u} - \frac{1}{u^{M+3}},$$

for $u$ large enough. Hence

(27) $$u \leq \frac{1}{1 - F_t(u)} \leq \frac{u}{1 - u^{-(M+2)}},$$

and for large $v$,

$$v \geq U_t(v) \geq v - v^{-M}.$$

Now we consider (19) and note that $\frac{x^0 - 1}{0}$ should be read as usual as $\log x$. We see that, with $a_t/U_t \equiv 1$,

$$\log U_t(vx) - \log U_t(v) - \log x$$
$$\leq \log vx - \log v - \log(1 - v^{-(M+1)}) - \log x \leq 2v^{-(M+1)}$$

and similarly,

$$-\log U_t(vx) + \log U_t(v) + \log x \leq 2(vx)^{-(M+1)}.$$

This implies (19) with $A_t(v) = v^{-M}$ and $\rho \equiv -\infty$.

Finally, we have to show (16), which has to be proved for the transformed process $\zeta_1$. As we see from (27), for large values of $v$ this transformation is very close to the identity function. So the transformed and untransformed processes are very close for high values. Nevertheless, the proof of (16) for the transformed process is more cumbersome than that for the untransformed process. We therefore confine ourselves to proving (16) for the untransformed process, since this proof contains the main ideas. Also, we will use the modified version of $E_{s,\delta}$, as described below Theorem 2.1, but we will keep the same notation. We have, for large enough $v$,

$$P\left\{\zeta_1 \notin E_{s,\delta} \Big| \sup_{t \in [s, s+\delta]} \zeta_1(t) \geq v\right\}$$
$$\leq P\left\{\zeta_1 \notin E_{s,\delta}, \sup_{t \in [s, s+\delta]} \zeta_1(t) \geq v\right\} \Big/ P\{\zeta_1(s) \geq v\}$$
$$\leq 2vP\left\{\zeta_1 \notin E_{s,\delta}, \sup_{t \in [s, s+\delta]} \zeta_1(t) \geq v\right\}$$

$$\leq 2vP\left\{\zeta_1 \notin E_{s,\delta}, \zeta_1(s) \geq \frac{v}{2}\right\} + 2vP\left\{\sup_{t\in[s,s+\delta]} \zeta_1(t) - \zeta_1(s) \geq \frac{v}{2}\right\}$$

$$=: D_1 + D_2.$$

Consider $D_1$ and use the independent increments property of a Wiener process:

$$D_1 = 2vP\left\{\sup_{t\in[s,s+\delta]} |\log(Ye^{W(t)-t/2}) - \log(Ye^{W(s)-s/2})| > K\left(\log\frac{1}{\delta}\right)^{-3},\right.$$

$$\left. Ye^{W(s)-s/2} \geq \frac{v}{2}\right\}$$

$$= 2vP\left\{\sup_{t\in[s,s+\delta]} |W(t) - W(s) - (t/2 - s/2)| > K\left(\log\frac{1}{\delta}\right)^{-3}\right\}$$

$$\times P\left\{Ye^{W(s)-s/2} \geq \frac{v}{2}\right\}$$

$$\leq 4P\left\{\sup_{t\in[s,s+\delta]} |W(t) - W(s)| > \frac{K}{2}\left(\log\frac{1}{\delta}\right)^{-3}\right\} \leq \exp\left(-\frac{1}{\sqrt{\delta}}\right),$$

where, for the last inequality, one of the well-known bounds for the oscillations of the Wiener process is used. For $D_2$ we obtain again, by the independent increments property,

$$D_2 = 2vP\left\{Y \sup_{t\in[s,s+\delta]} (e^{W(t)-t/2} - e^{W(s)-s/2}) \geq \frac{v}{2}\right\}$$

$$\leq 2vP\left\{Ye^{W(s)-s/2} \sup_{t\in[s,s+\delta]} (e^{W(t)-t/2-(W(s)-s/2)} - 1) \geq \frac{v}{2}\right\}$$

$$\leq 2vP\left\{Ye^{W(s)-s/2} \sup_{y\in[0,1]} (e^{\sqrt{\delta}V(y)} - 1) \geq \frac{v}{2}\right\} =: D_3,$$

where $V$ is a standard Wiener process independent of $W$ and $Y$. So the three terms in the latter probability are independent. Recall that $Ee^{W(s)-s/2} = 1$. Now

$$E \sup_{y\in[0,1]} (e^{\sqrt{\delta}V(y)} - 1) = E(e^{\sqrt{\delta}\sup_{y\in[0,1]} V(y)} - 1) = \int_0^\infty e^{\sqrt{\delta}x} 2\phi(x)\,dx - 1,$$

with $\phi$ the standard normal density. A straightforward calculation shows that the latter quantity is equal to $2e^{\delta/2}(1 - \Phi(-\sqrt{\delta})) - 1 \leq \sqrt{\delta}$. Hence

$$D_3 = 2vEP\left\{Y \geq \frac{v}{2e^{W(s)-s/2}\sup_{y\in[0,1]}(e^{\sqrt{\delta}V(y)} - 1)}\right|$$

$$e^{W(s)-s/2} \sup_{y\in[0,1]} (e^{\sqrt{\delta}V(y)} - 1) \Big\}$$

$$\leq 2v\frac{2}{v} E e^{W(s)-s/2} E \sup_{y\in[0,1]} (e^{\sqrt{\delta}V(y)} - 1) \leq 4 \cdot 1 \cdot \sqrt{\delta}.$$

So $D_1 + D_2 \leq 5\sqrt{\delta}$. This is much smaller than the bound required in (16). Hence, we have proved that condition.

It should be observed that, for both examples, condition (16) trivially remains satisfied if we transform the process $\xi_1$ by transformations of the marginals by increasing, continuous functions. So as long as these transformations yield a process that satisfies the other conditions (including that the transformed process is an element of $C[0,1]$), we have a new process for which Theorem 2.3 is valid. In this way we can obtain processes with many different, and nonconstant, extreme value index functions.

### 3. Proofs.

PROOF OF THEOREM 2.1. We only give a proof for the case $c = 1$; for general $c > 0$ the proof is similar. For any $\beta \in [0, \frac{1}{2})$ define

$$f_{t,x} = \mathbb{1}_{C_{t,x}} x^\beta,$$
$$\mathcal{F} = \{f_{t,x} : 0 \leq t \leq 1, x \geq 1\}.$$

Also, define the random measures

$$\mathcal{Z}_{n,i} = \frac{1}{\sqrt{k}} \delta_{\zeta_i k/n};$$

$\mathcal{Z}_{n,i}$ is a random function on $\mathcal{F}$ with

$$\mathcal{Z}_{n,i}(f_{t,x}) = \frac{1}{\sqrt{k}} \mathbb{1}_{\{\zeta_i(t)k/n \geq x\}} x^\beta.$$

Then

$$x^\beta w_n(t,x) = \sum_{i=1}^n (\mathcal{Z}_{n,i}(f_{t,x}) - E\mathcal{Z}_{n,i}(f_{t,x})).$$

First we are going to prove the tightness of $\{\sum_{i=1}^n (\mathcal{Z}_{n,i}(f) - E\mathcal{Z}_{n,i}(f)), f \in \mathcal{F}\}$. We need the following version of Theorem 2.11.9 in [13] (note that, indeed, the middle condition there is not needed here).

DEFINITION 3.1. For any $\varepsilon > 0$, the bracketing number $N_{[\cdot]}(\varepsilon, \mathcal{F}, L_2^n)$ is the minimal number of sets $N_\varepsilon$ in a partition $\mathcal{F} = \bigcup_{j=1}^{N_\varepsilon} \mathcal{F}_{\varepsilon j}$ of the index set into sets $\mathcal{F}_{\varepsilon j}$ independent of $n$ such that, for every partitioning set $\mathcal{F}_{\varepsilon j}$,

$$(28) \qquad \sum_{i=1}^{n} E^* \sup_{f,g \in \mathcal{F}_{\varepsilon j}} |\mathcal{Z}_{n,i}(f) - \mathcal{Z}_{n,i}(g)|^2 \leq \varepsilon^2.$$

THEOREM 3.2. *For each $n$, let $\mathcal{Z}_{n,1}, \mathcal{Z}_{n,2}, \ldots, \mathcal{Z}_{n,n}$ be independent stochastic processes with finite second moments indexed by a totally bounded semimetric space $(\mathcal{F}, d)$. Suppose*

$$(29) \qquad \sum_{i=1}^{n} E^* \|\mathcal{Z}_{n,i}\|_{\mathcal{F}} \mathbb{1}_{\{\|\mathcal{Z}_{n,i}\|_{\mathcal{F}} > \lambda\}} \to 0 \qquad \text{for every } \lambda > 0,$$

*where $\|\mathcal{Z}_{n,i}\|_{\mathcal{F}} = \sup_{f \in \mathcal{F}} |\mathcal{Z}_{n,i}(f)|$ and*

$$(30) \qquad \int_0^{\delta_n} \sqrt{\log N_{[\cdot]}(\varepsilon, \mathcal{F}, L_2^n)} \, d\varepsilon \to 0 \qquad \text{for every } \delta_n \downarrow 0.$$

*Then the sequence $\sum_{i=1}^{n} (\mathcal{Z}_{n,i} - E\mathcal{Z}_{n,i})$ is asymptotically tight in $\ell^\infty(\mathcal{F})$ and converges weakly, provided the finite-dimensional distributions converge weakly.*

We can define $d$ on $\mathcal{F}$ by $d(f_{t,x}, f_{s,y}) = d((t,x), (s,y))$; see the first paragraph of Section 2. We first show briefly that our class of functions $\mathcal{F}$ is totally bounded under the metric $d$. We consider w.l.o.g. only the case $x \leq y$. Since $\nu$ is a finite, and hence, tight measure on $\{h \in C[0,1] : \sup_{t \in [0,1]} h(t) \geq 1\}$, we can, for any $\delta_1 > 0$, find a $\delta_2 > 0$ such that if $|t - s| \leq \delta_2$, then

$$\nu(C_{s,y} \setminus C_{t,x}) \leq \nu(C_{s,y} \setminus C_{s,y+\delta_1/2}) + \nu(C_{s,y+\delta_1/2} \setminus C_{t,x})$$
$$\leq \frac{1}{y} - \frac{1}{y + \delta_1/2} + \frac{1}{2}\delta_1 \leq \delta_1$$

and (hence), if $\frac{1}{x} - \frac{1}{y} \leq \delta_1$, then

$$\nu(C_{t,x} \setminus C_{s,y}) \leq \nu(C_{t,x} \setminus C_{t,y}) + \nu(C_{t,y} \setminus C_{s,y}) \leq \frac{1}{x} - \frac{1}{y} + \delta_1 \leq 2\delta_1.$$

Now we have, for $|t - s| \leq \delta_2$ and $\frac{1}{x} - \frac{1}{y} \leq \delta_1$,

$$d^2(f_{t,x}, f_{s,y})$$
$$= (y^\beta - x^\beta)^2 \nu(C_{t,x} \cap C_{s,y}) + x^{2\beta} \nu(C_{t,x} \setminus C_{s,y}) + y^{2\beta} \nu(C_{s,y} \setminus C_{t,x})$$
$$\leq (y^\beta - x^\beta)^2 \nu(C_{s,y}) + x^{2\beta}\left(\frac{1}{x} \wedge 2\delta_1\right) + y^{2\beta}\left(\frac{1}{y} \wedge \delta_1\right)$$
$$(31)$$
$$\leq \left(xy^\beta\left(\frac{1}{x} - \frac{1}{y}\right)\right)^2 \frac{1}{y} + x^{2\beta}\left(\frac{1}{x} \wedge 2\delta_1\right) + y^{2\beta}\left(\frac{1}{y} \wedge \delta_1\right)$$

$$\leq x^{1+2\beta}\left(\frac{1}{x}-\frac{1}{y}\right)^2 + x^{2\beta}\left(\frac{1}{x}\wedge 2\delta_1\right) + y^{2\beta}\left(\frac{1}{y}\wedge\delta_1\right)$$

$$\leq \delta_1^{1-2\beta} + 2\delta_1^{1-2\beta} + \delta_1^{1-2\beta} = 4\delta_1^{1-2\beta}.$$

So, since $1-2\beta > 0$, we see that for $\varepsilon > 0$ we can find a $\delta_1 > 0$ such that $d(f_{t,x}, f_{s,y}) \leq \varepsilon$, for $\frac{1}{x} - \frac{1}{y} \leq \delta_1$ and $|t-s| \leq \delta_2$. Since obviously $\mathcal{F}$ is totally bounded under the metric $d_0(f_{t,x}, f_{s,y}) = |\frac{1}{x} - \frac{1}{y}| + |t-s|$, the total boundedness under $d$ follows.

To prove (29), observe

$$\|\mathcal{Z}_{n,i}\|_{\mathcal{F}} = \frac{1}{\sqrt{k}}\sup_{0\leq t\leq 1}\left(\zeta_i(t)\frac{k}{n}\right)^{\beta}.$$

So

$$\sum_{i=1}^n E\|\mathcal{Z}_{n,i}\|_{\mathcal{F}}\mathbb{1}_{\{\|\mathcal{Z}_{n,i}\|_{\mathcal{F}}>\lambda\}}$$

$$= \frac{n}{\sqrt{k}}E\left(\sup_{0\leq t\leq 1}\zeta_i(t)\frac{k}{n}\right)^{\beta}\mathbb{1}_{\{\sup_{0\leq t\leq 1}\zeta_i(t)(k/n)>(\sqrt{k}\lambda)^{1/\beta}\}}$$

(32)
$$= \frac{n}{\sqrt{k}}\int_{(\sqrt{k}\lambda)^{1/\beta}}^{\infty} x^{\beta}\,dF_n(x)$$

$$= -\frac{n}{\sqrt{k}}x^{\beta}(1-F_n(x))\Big|_{(\sqrt{k}\lambda)^{1/\beta}}^{\infty}$$

$$+ \beta\frac{n}{\sqrt{k}}\int_{(\sqrt{k}\lambda)^{1/\beta}}^{\infty} x^{\beta-1}(1-F_n(x))\,dx,$$

where $1-F_n(x) = P\{\sup_{0\leq t\leq 1}\zeta_i(t)\frac{k}{n} \geq x\}$. Note that $P\{\sup_{0\leq t\leq 1}\zeta_i(t) \geq x\} = x^{-1}\nu_x(\{h \in C[0,1]: \sup_{0\leq t\leq 1} h(t) \geq 1\})$. Hence it follows from (4) that the function $x \mapsto P\{\sup_{0\leq t\leq 1}\zeta_i(t) \geq x\}$ is regularly varying at infinity with exponent $-1$, so

$$\lim_{u\to\infty}\frac{P\{\sup_{0\leq t\leq 1}\zeta_i(t) \geq ux\}}{P\{\sup_{0\leq t\leq 1}\zeta_i(t) \geq u\}} = \frac{1}{x}, \qquad x > 0.$$

Let $0 < \tau < 1$. Now it immediately follows from Potter's inequality (see, e.g., [2]) that, for large $n$ and $x \geq 1$,

$$\frac{(n/k)P\{\sup_{0\leq t\leq 1}\zeta_i(t) \geq (n/k)x\}}{(n/k)P\{\sup_{0\leq t\leq 1}\zeta_i(t) \geq (n/k)\}} \leq 2x^{\tau-1}.$$

Also, we have as $n \to \infty$,

$$\frac{n}{k}P\left\{\sup_{0\leq t\leq 1}\zeta_i(t)\frac{k}{n} \geq 1\right\} = \nu_{n/k}\left(\left\{f \in C[0,1]: \sup_{0\leq t\leq 1} f(t) \geq 1\right\}\right)$$

$$\to \nu\left(\left\{f \in C[0,1]: \sup_{0 \leq t \leq 1} f(t) \geq 1\right\}\right) =: \frac{C}{3},$$

for some positive, finite $C$. So for large $n$ and $x \geq 1$,

(33) $$1 - F_n(x) \leq C\frac{k}{n}x^{\tau-1}.$$

Hence, the right-hand side of (32) is bounded from above by

$$Ck^{(2\beta+\tau-1)/(2\beta)}\lambda^{(\beta+\tau-1)/\beta} + \beta C\sqrt{k}\int_{(\sqrt{k}\lambda)^{1/\beta}}^{\infty} x^{\beta+\tau-2}\,dx$$

$$= C\frac{1-\tau}{1-\beta-\tau}\lambda^{(\beta+\tau-1)/\beta}k^{(2\beta+\tau-1)/(2\beta)} \to 0,$$

for $\tau$ small enough, since $\beta < \frac{1}{2}$. That is (29).

Next we will prove (30). For any (small) $\varepsilon > 0$, let $a = \varepsilon^{3/(2\beta-1)}$, $\delta = \exp\{-\varepsilon^{-1}\}$ and $\theta = 1/(1 - K\varepsilon^3)$. Define

$$\mathcal{F}(a) = \{f_{t,x} \in \mathcal{F}, x > a\},$$
$$\mathcal{F}(l,j) = \{f_{t,x} \in \mathcal{F}, l\delta \leq t < (l+1)\delta, \theta^j \leq x \leq \theta^{j+1}\}.$$

Then we have the "partition" $\mathcal{F} = \mathcal{F}(a) \cup \bigcup_{l=0}^{\lceil 1/\delta \rceil} \bigcup_{j=0}^{\lfloor \log a/\log \theta \rfloor} \mathcal{F}(l,j)$. First we check (28) for $\mathcal{F}(a)$:

$$\sum_{i=1}^{n} E \sup_{f,g \in \mathcal{F}(a)} |\mathcal{Z}_{n,i}(f) - \mathcal{Z}_{n,i}(g)|^2$$

$$= nE \sup_{f,g \in \mathcal{F}(a)} (\mathcal{Z}_{n,i}(f) - \mathcal{Z}_{n,i}(g))^2$$

$$\leq 4nE \sup_{f \in \mathcal{F}(a)} \mathcal{Z}_{n,i}^2(f)$$

$$\leq \frac{4n}{k} E\left(\sup_{0 \leq t \leq 1} \zeta_i(t)\frac{n}{k}\right)^{2\beta} \mathbb{1}_{\{\sup_{0 \leq t \leq 1} \zeta_i(t)k/n \geq a\}}$$

$$= \frac{4n}{k}\int_{a}^{\infty} x^{2\beta}\,dF_n(x)$$

$$\leq 4C\frac{1-\tau}{1-2\beta-\tau}a^{2\beta+\tau-1}$$

$$= 4C\frac{1-\tau}{1-2\beta-\tau}\varepsilon^{3(2\beta+\tau-1)/(2\beta-1)},$$

where the last inequality follows from integration by parts and (33). Clearly, the latter expression is bounded from above by $\varepsilon^2$ for $\tau$ (and $\varepsilon$) small enough.

Now we consider (28) for the $\mathcal{F}(l,j)$. First note that

$$\sup_{f \in \mathcal{F}(l,j)} \mathcal{Z}_{n,i}(f) \leq \frac{1}{\sqrt{k}} \mathbb{1}_{\{\sup_{l\delta \leq t < (l+1)\delta} \zeta_i(t)(k/n) \geq \theta^j\}} \theta^{(j+1)\beta}$$

$$= \frac{1}{\sqrt{k}} \mathbb{1}_{\{\sup_{l\delta \leq t < (l+1)\delta} \zeta_i(t)(k/n) \geq \theta^j, \zeta_i \in E_{s,\delta}\}} \theta^{(j+1)\beta}$$

$$+ \frac{1}{\sqrt{k}} \mathbb{1}_{\{\sup_{l\delta \leq t < (l+1)\delta} \zeta_i(t)(k/n) \geq \theta^j, \zeta_i \notin E_{s,\delta}\}} \theta^{(j+1)\beta}.$$

Suppose $\zeta_i \in E_{s,\delta}$ and $\sup_{l\delta \leq t < (l+1)\delta} \zeta_i(t)\frac{k}{n} \geq \theta^j$. Then for small enough $\delta$,

$$\sup_{l\delta \leq t < (l+1)\delta} \zeta_i(t) - \zeta_i(l\delta) \leq K\varepsilon^3 \zeta_i(l\delta),$$

and, hence, $\zeta_i(l\delta)\frac{k}{n} \geq \theta^{j-1}$. So

$$\sup_{f \in \mathcal{F}(l,j)} \mathcal{Z}_{n,i}(f) \leq \frac{1}{\sqrt{k}} \mathbb{1}_{\{\zeta_i(l\delta)(k/n) \geq \theta^{j-1}, \zeta_i \in E_{s,\delta}\}} \theta^{(j+1)\beta}$$

$$+ \frac{1}{\sqrt{k}} \mathbb{1}_{\{\sup_{l\delta \leq t < (l+1)\delta} \zeta_i(t)(k/n) \geq \theta^j, \zeta_i \notin E_{s,\delta}\}} \theta^{(j+1)\beta}.$$

Similarly, it can be shown that

$$\inf_{f \in \mathcal{F}(l,j)} \mathcal{Z}_{n,i}(f) \geq \frac{1}{\sqrt{k}} \mathbb{1}_{\{\zeta_i(l\delta)(k/n) \geq \theta^{j+2}, \zeta_i \in E_{s,\delta}\}} \theta^{j\beta}.$$

This yields

$$\sum_{i=1}^{n} E^* \sup_{f,g \in \mathcal{F}(l,j)} |\mathcal{Z}_{n,i}(f) - \mathcal{Z}_{n,i}(g)|^2$$

$$\leq nE^* \left( \sup_{f \in \mathcal{F}(l,j)} \mathcal{Z}_{n,i}(f) - \inf_{f \in \mathcal{F}(l,j)} \mathcal{Z}_{n,i}(f) \right)^2$$

$$\leq \frac{n}{k} E(\mathbb{1}_{\{\zeta_i(l\delta)(k/n) \geq \theta^{j-1}, \zeta_i \in E_{s,\delta}\}} \theta^{(j+1)\beta}$$

$$+ \mathbb{1}_{\{\sup_{l\delta \leq t < (l+1)\delta} \zeta_i(t)(k/n) \geq \theta^j, \zeta_i \notin E_{s,\delta}\}} \theta^{(j+1)\beta}$$

$$- \mathbb{1}_{\{\zeta_i(l\delta)(k/n) \geq \theta^{j+2}, \zeta_i \in E_{s,\delta}\}} \theta^{j\beta})^2$$

$$= \frac{n}{k} E(\mathbb{1}_{\{\zeta_i(l\delta)(k/n) \geq \theta^{j-1}, \zeta_i \in E_{s,\delta}\}} \theta^{(j+1)\beta} - \mathbb{1}_{\{\zeta_i(l\delta)(k/n) \geq \theta^{j+2}, \zeta_i \in E_{s,\delta}\}} \theta^{j\beta})^2$$

$$+ \frac{n}{k} P \left\{ \sup_{l\delta \leq t < (l+1)\delta} \zeta_i(t) \frac{k}{n} \geq \theta^j, \zeta_i \notin E_{s,\delta} \right\} \theta^{2(j+1)\beta}$$

$$=: T_1 + T_2.$$

We have

$$T_1 \leq \frac{n}{k} E((\theta^{(j+1)\beta} - \theta^{j\beta})\mathbb{1}_{\{\zeta_i(l\delta)(k/n) \geq \theta^{j-1}\}} + \theta^{j\beta}\mathbb{1}_{\{\theta^{j+2} > \zeta_i(l\delta)(k/n) \geq \theta^{j-1}\}})^2$$

$$\leq 2\theta^{2(j+1)\beta}(1-\theta^{-\beta})^2 \frac{1}{\theta^{j-1}} + 2\theta^{2j\beta}\left(\frac{1}{\theta^{j-1}} - \frac{1}{\theta^{j+2}}\right)$$

$$\leq 2\theta^{j+1}(1-\theta^{-1/2})^2 \frac{1}{\theta^{j-1}} + 2\theta^j\left(\frac{1}{\theta^{j-1}} - \frac{1}{\theta^{j+2}}\right)$$

$$\leq 3(K\varepsilon^3 + 3K\varepsilon^3) \leq \frac{1}{2}\varepsilon^2,$$

and for large $n$,

$$T_2 \leq \frac{n}{k} P\left\{\sup_{t \in [l\delta,(l+1)\delta)} \zeta_i(t)\frac{k}{n} \geq 1, \zeta_i \notin E_{s,\delta}\right\}\theta^{2(j+1)\beta}$$

$$\leq \frac{n}{k} P\left\{\sup_{t \in [l\delta,(l+1)\delta)} \zeta_i(t)\frac{k}{n} \geq 1\right\} P\left\{\zeta_i \notin E_{s,\delta}\Big| \sup_{t \in [l\delta,(l+1)\delta)} \zeta_i(t) \geq \frac{n}{k}\right\} a^{2\beta}\theta^{2\beta}$$

$$\leq C \cdot c_1 \left(\log \frac{1}{\delta}\right)^{-(2+2\beta)/(1-2\beta)} a^{2\beta} \leq \frac{1}{2}\varepsilon^2.$$

Hence, we have shown (28).

It is easy to see that the number of elements of the partition is bounded by $\exp(2/\varepsilon)$, which leads to (30). Hence, by Theorem 3.2 we have proved the asymptotic tightness condition.

It remains to prove that the finite-dimensional distributions of $\sum_{i=1}^{n}(\mathcal{Z}_{n,i} - E\mathcal{Z}_{n,i})$ converge weakly. This follows from the fact that multivariate weak convergence follows from weak convergence of linear combinations of the components and the (univariate) Lindeberg–Feller central limit theorem. It is easily seen that the Lindeberg condition is fulfilled for the linear combinations, since the $f_{t,x}$ are made up of indicators and hence bounded.

The fact that $Z$ is bounded and uniformly $d$-continuous follows from the general theory of weak convergence and properties of Gaussian processes; see Section 1.5 in [13]. $\square$

PROOF OF COROLLARY 2.2. Write $V_{n,t} = \zeta_{n-k,n}(t)\frac{k}{n}$. We first show the result for $\alpha = -1$, that is,

(34) $$\sup_{0 \leq t \leq 1} \left|\sqrt{k}\left(\frac{1}{V_{n,t}} - 1\right) + W(C_{t,1})\right| \xrightarrow{P} 0.$$

Clearly,

$$\sup_{0 \leq t \leq 1} \left|\sqrt{k}\left(\frac{1}{V_{n,t}} - 1\right) + w_n(t, V_{n,t})\right| \xrightarrow{P} 0,$$

so (17), with $\beta = 0$, yields

$$\sup_{0 \leq t \leq 1} \left| \sqrt{k}\left(\frac{1}{V_{n,t}} - 1\right) + W(C_{t,V_{n,t}}) \right| \xrightarrow{P} 0.$$

Now by the boundedness and uniform $d$-continuity of $W$ we obtain (34). Finally, write

$$\sqrt{k}(V_{n,t}^{\alpha(t)} - 1) = \sqrt{k}(V_{n,t}^{-1} - 1) \frac{V_{n,t}^{\alpha(t)} - 1}{V_{n,t}^{-1} - 1}.$$

Since, by (34),

$$\sup_{0 \leq t \leq 1} \left| \frac{V_{n,t}^{\alpha(t)} - 1}{V_{n,t}^{-1} - 1} + \alpha(t) \right| \xrightarrow{P} 0,$$

we obtain, again using (34), (18). □

PROOF OF THEOREM 2.3. First, from (19) we can prove, for any $\varepsilon > 0$, there exists $s_\varepsilon > 0$ such that if $v > v_\varepsilon$ and $x \geq 1$, we have, for all $0 \leq t \leq 1$,

$$(35) \quad \left| \left( \frac{\log U_t(vx) - \log U_t(v)}{a_t(v)/U_t(v)} - \frac{x^{\gamma^-(t)} - 1}{\gamma^-(t)} \right) \bigg/ A_t(v) - H_{\gamma^-(t), \rho(t)}(x) \right|$$
$$\leq \varepsilon(1 + x^{\gamma^-(t) + \varepsilon});$$

the proof follows along the lines of that for the one-dimensional situation in [3]; see also [8]. Inequality (35) implies

$$(36) \quad \left| \frac{\log U_t(vx) - \log U_t(v)}{a_t(v)/U_t(v)} - \frac{x^{\gamma^-(t)} - 1}{\gamma^-(t)} \right| \leq |A_t(v)|(C_\varepsilon + x^\varepsilon),$$

where $C_\varepsilon \in (0, \infty)$ is a constant. Note that

$$M_n^{(1)}(t) = \frac{1}{k} \sum_{i=1}^{k-1} \log U_t(\zeta_{n-i,n}(t)) - \log U_t(\zeta_{n-k,n}(t)).$$

Hence, we have, for sufficiently large $n$,

$$\frac{1}{k} \sum_{i=0}^{k-1} \frac{(\zeta_{n-i,n}(t)/\zeta_{n-k,n}(t))^{\gamma^-(t)} - 1}{\gamma^-(t)}$$
$$- |A_t(\zeta_{n-k,n}(t))| \left( C_\varepsilon + \frac{1}{k} \sum_{i=0}^{k-1} \left( \frac{\zeta_{n-i,n}(t)}{\zeta_{n-k,n}(t)} \right)^\varepsilon \right)$$
$$(37) \qquad \leq \frac{M_n^{(1)}(t)}{a_t(\zeta_{n-k,n}(t))/U_t(\zeta_{n-k,n}(t))}$$

$$\leq \frac{1}{k}\sum_{i=0}^{k-1} \frac{(\zeta_{n-i,n}(t)/\zeta_{n-k,n}(t))^{\gamma^-(t)} - 1}{\gamma^-(t)}$$

$$+ |A_t(\zeta_{n-k,n}(t))|\left(C_\varepsilon + \frac{1}{k}\sum_{i=0}^{k-1}\left(\frac{\zeta_{n-i,n}(t)}{\zeta_{n-k,n}(t)}\right)^\varepsilon\right).$$

As before, write $V_{n,t} = \zeta_{n-k,n}(t)\frac{k}{n}$. Next

$$\sqrt{k}\left(\frac{1}{k}\sum_{i=0}^{k-1}\frac{(\zeta_{n-i,n}(t)/\zeta_{n-k,n}(t))^{\gamma^-(t)} - 1}{\gamma^-(t)} - \frac{1}{1-\gamma^-(t)}\right)$$

$$= \sqrt{k}\left(\int_{V_{n,t}}^\infty \frac{(x/V_{n,t})^{\gamma^-(t)} - 1}{\gamma^-(t)}\, d\left(-\frac{n}{k}S_{n,t}\left(x\frac{n}{k}\right)\right) - \frac{1}{1-\gamma^-(t)}\right)$$

$$= \sqrt{k}\left(V_{n,t}^{-\gamma^-(t)}\int_{V_{n,t}}^\infty \frac{n}{k}S_{n,t}\left(x\frac{n}{k}\right)x^{\gamma^-(t)-1}\,dx - \int_1^\infty x^{\gamma^-(t)-2}\,dx\right)$$

$$= V_{n,t}^{-\gamma^-(t)}\int_{V_{n,t}}^\infty w_n(t,x)x^{\gamma^-(t)-1}\,dx$$

$$+ \sqrt{k}(V_{n,t}^{-\gamma^-(t)} - 1)\int_{V_{n,t}}^\infty x^{\gamma^-(t)-2}\,dx + \sqrt{k}\int_{V_{n,t}}^1 x^{\gamma^-(t)-2}\,dx.$$

So

$$\sqrt{k}\left(\frac{1}{k}\sum_{i=0}^{k-1}\frac{(\zeta_{n-i,n}(t)/\zeta_{n-k,n}(t))^{\gamma^-(t)} - 1}{\gamma^-(t)} - \frac{1}{1-\gamma^-(t)}\right) - \mathcal{P}(t)$$

$$= V_{n,t}^{-\gamma^-(t)}\int_{V_{n,t}}^\infty (w_n(t,x) - W(C_{t,x}))x^{\gamma^-(t)-1}\,dx$$

$$+ (V_{n,t}^{-\gamma^-(t)} - 1)\int_{Vn,t}^\infty W(C_{t,x})x^{\gamma^-(t)-1}\,dx$$

(38)
$$+ (\sqrt{k}(V_{n,t}^{-\gamma^-(t)} - 1) + \gamma^-(t)W(C_{t,1}))\int_{V_{n,t}}^\infty x^{\gamma^-(t)-2}\,dx$$

$$+ \left(\sqrt{k}\int_{V_{n,t}}^1 x^{\gamma^-(t)-2}\,dx + W(C_{t,1})\right)$$

$$- \int_1^{V_{n,t}} W(C_{t,x})x^{\gamma^-(t)-1}\,dx + \gamma^-(t)W(C_{t,1})\int_1^{V_{n,t}} x^{\gamma^-(t)-2}\,dx.$$

From Theorem 2.1 we obtain, for the first term on the right-hand side in (38),

$$\sup_{t\in[0,1]} V_{n,t}^{-\gamma^-(t)}\left|\int_{V_{n,t}}^\infty (w_n(t,x) - W(C_{t,x}))x^{\gamma^-(t)-1}\,dx\right|$$

$$
\begin{aligned}
(39) \quad &\leq \sup_{t\in[0,1]} V_{n,t}^{-\gamma^-(t)} \cdot \sup_{t\in[0,1], x\geq V_{n,t}} x^\beta |w_n(t,x) - W(C_{t,x})| \\
&\quad \times \sup_{t\in[0,1]} \int_{V_{n,t}}^\infty y^{\gamma^-(t)-1-\beta} \, dy.
\end{aligned}
$$

Now it follows from Theorem 2.1 with $\beta$ *positive* (this is crucial) and Corollary 2.2 that the right-hand side of (39) converges to 0 in probability. It readily follows from Corollary 2.2 that the five other terms on the right-hand side of (38) converge to 0 in probability. So we have

$$
(40) \quad \sup_{0\leq t\leq 1}\left|\sqrt{k}\left(\frac{1}{k}\sum_{i=0}^{k-1}\frac{(\zeta_{n-i,n}(t)/\zeta_{n-k,n}(t))^{\gamma^-(t)} - 1}{\gamma^-(t)} - \frac{1}{1-\gamma^-(t)}\right) - \mathcal{P}(t)\right| \xrightarrow{P} 0
$$

as $n\to\infty$.

For the remainder term of

$$
\frac{M_n^{(1)}(t)}{a_t(\zeta_{n-k,n}(t))/U_t(\zeta_{n-k,n}(t))}
$$

in (37), note that we obtain from Lemma 3.2 in [5] that, for $0\leq \varepsilon < 1$,

$$
(41) \quad \sup_{0\leq t\leq 1}\left|\frac{1}{k}\sum_{i=0}^{k-1}\left(\frac{\zeta_{n-i,n}(t)}{\zeta_{n-k,n}(t)}\right)^\varepsilon - \frac{1}{1-\varepsilon}\right| \xrightarrow{P} 0 \qquad \text{as } n\to\infty.
$$

It can be derived from the second-order condition (19) and Corollary 2.2 that

$$
\sup_{0\leq t\leq 1}\left|\frac{A_t(n/k)}{A_t(\zeta_{n-k,n}(t))} - 1\right| \xrightarrow{P} 0.
$$

Using this in combination with (20) and (41), we see that the remainder term in (37) is negligible, so we obtain that

$$
(42) \quad \sup_{0\leq t\leq 1}\left|\sqrt{k}\left(\frac{M_n^{(1)}(t)}{a_t(\zeta_{n-k,n}(t))/U_t(\zeta_{n-k,n}(t))} - \frac{1}{1-\gamma^-(t)}\right) - \mathcal{P}(t)\right| \xrightarrow{P} 0
$$

as $n\to\infty$. Similarly,

$$
(43) \quad \sup_{0\leq t\leq 1}\left|\sqrt{k}\left(\frac{M_n^{(2)}(t)}{(a_t(\zeta_{n-k,n}(t))/U_t(\zeta_{n-k,n}(t)))^2} - \frac{2}{(1-\gamma^-(t))(1-2\gamma^-(t))}\right) - \mathcal{Q}(t)\right| \xrightarrow{P} 0
$$

as $n \to \infty$. Hence, we get

(44) $$\sup_{0 \le t \le 1} |\sqrt{k}(\hat{\gamma}_n^-(t) - \gamma^-(t)) - \mathcal{M}(t)| \xrightarrow{P} 0$$

as $n \to \infty$, where

$$\mathcal{M}(t) = -2(1-\gamma^-(t))^2(1-2\gamma^-(t))\mathcal{P}(t) + \tfrac{1}{2}(1-\gamma^-(t))^2(1-2\gamma^-(t))^2\mathcal{Q}(t).$$

We now prove (22). Write

$$\sqrt{k}(\hat{\gamma}_n^+(t) - \gamma^+(t)) - \gamma^+(t)\mathcal{P}(t)$$

$$= \frac{a_t(\zeta_{n-k,n}(t))}{U_t(\zeta_{n-k,n}(t))}\left(\sqrt{k}\left(\frac{M_n^{(1)}(t)}{a_t(\zeta_{n-k,n}(t))/U_t(\zeta_{n-k,n}(t))} - \frac{1}{1-\gamma^-(t)}\right) - \mathcal{P}(t)\right)$$

$$+ \sqrt{k}\left(\frac{a_t(\zeta_{n-k,n}(t))}{U_t(\zeta_{n-k,n}(t))} - \gamma^+(t)\right)\frac{1}{1-\gamma^-(t)}$$

$$+ \left(\frac{a_t(\zeta_{n-k,n}(t))}{U_t(\zeta_{n-k,n}(t))} - \gamma^+(t)\right)\mathcal{P}(t).$$

If we show that

(45) $$\sup_{0 \le t \le 1} \left|\sqrt{k}\left(\frac{a_t(\zeta_{n-k,n}(t))}{U_t(\zeta_{n-k,n}(t))} - \gamma^+(t)\right)\right| \xrightarrow{P} 0,$$

then (22) follows from (42). We have

$$\sqrt{k}\left(\frac{a_t(\zeta_{n-k,n}(t))}{U_t(\zeta_{n-k,n}(t))} - \gamma^+(t)\right)$$

$$= \sqrt{k}\left(\frac{a_t(n/k)}{U_t(n/k)} - \gamma^+(t)\right)\frac{a_t(\zeta_{n-k,n}(t))/U_t(\zeta_{n-k,n}(t))}{a_t(n/k)/U_t(n/k)}$$

$$+ \sqrt{k}\left(\frac{a_t(\zeta_{n-k,n}(t))/U_t(\zeta_{n-k,n}(t))}{a_t(n/k)/U_t(n/k)} - \left(\zeta_{n-k,n}(t)\frac{k}{n}\right)^{\gamma^-(t)}\right)\gamma^+(t)$$

$$+ \sqrt{k}\left(\left(\zeta_{n-k,n}(t)\frac{k}{n}\right)^{\gamma^-(t)} - 1\right)\gamma^+(t).$$

From [8] and [10] it follows that (19) implies

(46) $$\frac{(a_t(xv)/U_t(xv))/(a_t(v)/U_t(v)) - x^{\gamma^-(t)}}{A_t(v)} \to x^{\gamma^-(t)}\frac{x^{\rho(t)} - 1}{\rho(t)}$$

as $v \to \infty$,

uniformly in $t \in [0,1]$ and locally uniformly in $x > 0$. Using (21), (20) and Corollary 2.2, we indeed obtain (45) and, hence, we have proved (22). Finally, we obtain (23) from (22) and (44).

For (24), note

$$\sqrt{k}\frac{\hat{U}_t(n/k) - U_t(n/k)}{a_t(n/k)}$$
$$= \sqrt{k}\frac{\log U_t(\zeta_{n-k,n}(t)) - \log U_t(n/k)}{a_t(n/k)/U_t(n/k)}$$
$$\times \left(\log\left(\frac{\xi_{n-k,n}(t)}{U_t(n/k)}\right)\right)^{-1}\left(\frac{\xi_{n-k,n}(t)}{U_t(n/k)} - 1\right)$$

and

$$\frac{\xi_{n-k,n}(t)}{U_t(n/k)} - 1 = \frac{\xi_{n-k,n}(t) - U_t(n/k)}{a_t(n/k)}\frac{a_t(n/k)}{U_t(n/k)}.$$

From Lemma 3.4 in [5] we obtain

$$\lim_{n\to\infty}\sup_{0\le t\le 1}\left|\frac{a_t(n/k)}{U_t(n/k)} - \gamma^+(t)\right| = 0.$$

Combining this with (14) yields

$$\sup_{0\le t\le 1}\left|\frac{\xi_{n-k,n}(t)}{U_t(n/k)} - 1\right| \xrightarrow{P} 0 \qquad \text{as } n\to\infty.$$

Hence,

$$\sup_{0\le t\le 1}\left|\left(\log\left(\frac{\xi_{n-k,n}(t)}{U_t(n/k)}\right)\right)^{-1}\left(\frac{\xi_{n-k,n}(t)}{U_t(n/k)} - 1\right) - 1\right| \xrightarrow{P} 0 \qquad \text{as } n\to\infty.$$

A proof similar to the one leading to (42) shows

$$\sup_{0\le t\le 1}\left|\sqrt{k}\frac{\log U_t(\zeta_{n-k,n}(t)) - \log U_t(n/k)}{a_t(n/k)/U_t(n/k)} - \mathcal{U}(t)\right| \xrightarrow{P} 0 \qquad \text{as } n\to\infty.$$

So we have obtained (24).

For (25), we use

$$\sqrt{k}\left(\frac{\hat{a}_t(n/k)}{a_t(n/k)} - 1\right)$$
$$= \sqrt{k}\left(\frac{M_n^{(1)}}{a_t(\zeta_{n-k,n})/U_t(\zeta_{n-k,n})} - \frac{1}{1-\gamma^-(t)}\right)(1-\hat{\gamma}_n^-(t))\frac{a_t(\zeta_{n-k,n}(t))}{a_t(n/k)}$$
$$- \sqrt{k}(\hat{\gamma}_n^-(t) - \gamma^-(t))\frac{1}{1-\gamma^-(t)}\frac{a_t(\zeta_{n-k,n}(t))}{a_t(n/k)}$$
$$+ \sqrt{k}\left(\frac{a_t(\zeta_{n-k,n}(t))}{a_t(n/k)} - 1\right).$$

Now

$$\sqrt{k}\bigg(\frac{a_t(\zeta_{n-k,n}(t))}{a_t(n/k)}-1\bigg)$$

$$=\sqrt{k}\bigg(\frac{a_t(\zeta_{n-k,n}(t))/U_t(\zeta_{n-k,n}(t))}{a_t(n/k)/U_t(n/k)}-\bigg(\frac{k}{n}\zeta_{n-k,n}(t)\bigg)^{\gamma^-(t)}\bigg)\frac{U_t(\zeta_{n-k,n}(t))}{U_t(n/k)}$$

$$+\sqrt{k}\bigg(\bigg(\frac{k}{n}\zeta_{n-k,n}(t)\bigg)^{\gamma^-(t)}-1\bigg)\frac{U_t(\zeta_{n-k,n}(t))}{U_t(n/k)}$$

$$+\sqrt{k}\frac{U_t(\zeta_{n-k,n}(t))-U_t(n/k)}{a_t(n/k)}\frac{a_t(n/k)}{U_t(n/k)}.$$

From (46) and (20), we know the first term tends to 0 in probability, uniformly in $t \in [0,1]$. Hence Corollary 2.2 and (24) yield

$$\sup_{0\leq t\leq 1}\bigg|\sqrt{k}\bigg(\frac{a_t(\zeta_{n-k,n}(t))}{a_t(n/k)}-1\bigg)-\gamma(t)W(C_{t,1})\bigg|\xrightarrow{P}0.$$

Using (42), (44) and Theorem 1.1, (25) now follows. □

**Acknowledgments.** The research of the second named author was performed at Erasmus University, Rotterdam, and EURANDOM, Eindhoven. We are grateful to Laurens de Haan for stimulating interest during the preparation of the paper and for several useful comments. We also thank two referees and the Coeditor for their constructive remarks.

Department of Econometrics
and Operations Research
Tilburg University
P.O. Box 90153
5000 LE Tilburg
The Netherlands
E-mail: j.h.j.einmahl@uvt.nl

Department of Mathematics
Xiamen University
Postcode 361005
Xiamen
China
E-mail: weiminghult@tom.com